\documentclass{article}
\usepackage{amssymb}
\title{On the Cohomology of spatial polygons in Euclidean spaces}
\newcommand{\ex}{\textbf{Example :}}
\newcommand{\nn}{\nonumber}
\newenvironment{defn}{\textbf{Definition :}}

\author{Vehbi Emrah Paksoy}
\date{}
\begin{document}
\maketitle
\begin{abstract}
Space of spatial polygons in Euclidean spaces has been studied extensively in \cite{KLY, HK, KM1}. There is a beautiful description of cohomology in \cite{HK}. In this paper we introduce another, easy to compute method to obtain the cohomology without using toric variety arguments. We also give a criterion for a polygon space to be Fano, thus, having ample anticanonical class. \footnote {I want to thank A. Klyachko for his ideas and the wonderful course "Geometric Invariant Theory" he taught at Bilkent University, Turkey.}
\end{abstract}
\section{Introduction and definitions}
Let $\mathcal{P}_{n}$ be the space of all $n$-gons with distinguished vertices in Euclidean space $\mathbb{E}^{3}$. An $n$-gon $P$ is determined by its vertices $v_{1},\ldots,v_{n}$. These vertices are joined in cyclic order by edges $p_{1},\ldots,p_{n}$ where $p_{i}$ is the oriented line segment from $v_{i}$ to $v_{i+1}$. Two polygons $P,Q \in \mathcal{P}_{n}$ are identified if and only if there exists an orientation preserving isometry $g$ of $\mathbb{E}^{3}$ which sends the vertices of $P$ to the vertices of $Q$.\\
\begin{defn} Let $m=(m_{1},\ldots,m_{n})$ be an $n$-tuple of positive real numbers. Then $\mathcal{M}_{n}$(or just $\mathcal{M}$) is the space of $n$-gons with side lenghts $m_{1},\ldots,m_{n}$ modulo isometries as above.
\end{defn}
Note that for any $P \in \mathcal{P}_{n}$, the vector of lenghts of sides of $P$ satisfies the following;
$$m_{i}<m_{1}+m_{2}+\cdots +\widehat{m_{i}}+\cdots + m_{n};\qquad i=\overline{1,n}=1,\ldots ,n.$$

It is known that $\mathcal{M}$ has only isolated singularities corresponding to the degenerate polygons. It is non singular if all sums $\pm m_{1} \pm m_{2}\pm \cdots \pm m_{n}$ are non zero. The following theorem about the structure of $\mathcal{M}$ appears in \cite{KLY}.\\
\textbf{Theorem :} If lenghts of all sides $\Vert p_{i} \Vert =m_{i}$ are rational and  $m_{1} \pm m_{2} \pm\cdots \pm m_{n} \ne 0$ then $\mathcal{M}$ is a non-singular projective variety.

We are going to use the relation between spatial polygons and stable weighted point configurations on the complex projective line. We have the following description;\\
\begin{defn} A $n$-point configuration $\Sigma$ is a collection of $n$-points $p_{1},\ldots,p_{n} \in \mathbb{P}^{1}$. Assume there is a given positive weight $m_{i}$ for each point.
The configuration of weighted points is called semi-stable (resp. stable) if sum of the weights of equal points does not exceed(resp. less than) half the weight of all points.
\end{defn}

Using Hilbert-Mumford stability criterion(\cite{M}), we can say that there exist a non-singular geometric factor of stable configurations with respect to a natural action of $PSL_{2}(\mathbb{C})$. It will be denoted by $\mathcal{C}_{n}(m)$ where  $m=(m_{1},\ldots,m_{n})$ is the vector of weights. By definition, $\mathcal{C}_{n}(m)$ is non-empty iff the weights satisfy the following polygon inequality
$$m_{i}<m_{1}+m_{2}+\ldots +\widehat{m_{i}}+\ldots+m_{n} ;\qquad i\in \{1,\ldots,n \} \qquad(*)$$\\
In a similar way, there exist a categorical factor of space of semi-stable configurations denoted by $\overline{\mathcal{C}}_{n}(m).$\\

Under condition (*), the variety $\overline{\mathcal{C}}_{n}(m)$ is a projective compactification of $\mathcal{C}_{n}(m)$ by a finite number of points. Its ample sheaf $\mathcal{O}(1)$ and the corresponding line bundle $\mathcal{L}$ may be described as follows. Let $\mathcal{T}(p_{i})$ be a tangent space at the point $p_{i} \in \mathbb{P}^{1}$. Then $\mathcal{L}$ is a line bundle on $\overline{\mathcal{C}}_{n}(m)$ with fiber 
$$\mathcal{L}(\mathbf{p})=\mathcal{T}(p_{1})^{\otimes m_{1}} \otimes \mathcal{T}(p_{2})^{\otimes m_{2}}\otimes \cdots \otimes \mathcal{T}(p_{n})^{\otimes m_{n}}$$
at a point $\mathbf{p}=(p_{1},\ldots,p_{n}) \in \overline{\mathcal{C}}_{n}(m).$
If all semi-stable configurations of weight $m$ are stable then $\overline{\mathcal{C}}_{n}(m)$=$\mathcal{C}_{n}(m)$ is a non-singular projective variety of dimension $n-3$.\\
\ex Let all weights $m_{i}=1$ i.e, $m=(1,\ldots,1)$. Then  $\overline{\mathcal{C}}_{n}(m)=\mathcal{C}_{n}(m)$ is a non-singular projective variety for odd $n$. In this case all sums $m_{1}\pm m_{2} \pm \cdots \pm m_{n}$ are non-zero.\\
\ex Let $\Sigma=(p_{1},\ldots,p_{n})$ be a configuration of $n$-points in $\mathbb{P}^{1}$ having one massive point, say $m_{1}$ i.e, $m_{1}+m_{i}>\frac{m}{2},m=m_{1}+\cdots +m_{n}$ so that $p_{1}\neq p_{i} ,\linebreak \forall i\neq 1$. Then we can interchange the coordinates in $\mathbb{P}^{1}$ such that $p_{1}=\infty,\linebreak p_{2}=0 \mbox{ and } p_{i}=z_{i},z_{i}\in \mathbb{C}; \quad z_{i} $ are defined uniquely up to multiplication by scalar multiplication $z\mapsto \lambda z$ which preserves $\infty,0$. Then moduli of the configuration is eqivalent to
$$\{(z_{3}:\ldots :z_{n})|z_{i}\in \mathbb{C} ,i=3,\ldots,n;\mbox{ not all zero}\}=\mathbb{P}^{n-3}$$\\
\ex Let $\Sigma=(p_{1},\ldots,p_{n})$ be a configuration with three massive points 
$$m_{i}+m_{j}>\frac{m}{2},m_{j}+m_{k}>\frac{m}{2},m_{i}+m_{k}>\frac{m}{2}$$\\
Then $p_{i}\neq p_{j},p_{j} \neq p_{k}$ and $p_{i}\neq p_{k}$. By a suitable coordinate change we may take $p_{i}=0,p_{j}=1,p_{k}=\infty$ and hence the moduli of configuration is equivalent to 
$$\prod_{\alpha \neq i,j,k} \mathbb{P}^{1}=(\mathbb{P}^{1})^{n-3}.$$

The following theorem reveals the relations between the variety of spatial polygons in $\mathbb{E}^{3}$ and stable configurations on $\mathbb{P}^{1}$.\\
\textbf{Theorem :}The algebraic variety of spatial polygons $\mathcal{M}$ is biregular equivalent to the variety $\overline{\mathcal{C}}_{n}(m)$ of semi-stable weighted configurations of points in the projective line.\\
{\bf Proof:} See \cite{KLY1} \qquad $\Box$\\
\section{Cohomology}
Algebraic variety of spatial polygons in Euclidean space $\mathbb{E}^{3}$ is equivalent to $\mathcal{C}_{n}(m)$, stable weighted configurations on complex projective line $\mathbb{P}^{1}=S^{2}$ modulo M\"obius group $PSL_{2}(\mathbb{C})$.\\
Let`s define $\mathcal{L}_{i}$ to be the linear vector bundle over $\mathcal{C}_{n}(m)$ such that fiber at $\Sigma=(p_{1},\ldots,p_{n})$ is equal to tangent space at $p_{i}\in \mathbb{P}^{1},i=1,\ldots,n.$ We call $\mathcal{L}_{i}$`s \textit{ natural bundles } on  $\mathcal{C}_{n}(m)$.\\
We know that  $\mathcal{C}_{n}(m)$ is non-empty if and only if weights $m=(m_{1},\ldots,m_{n})$ satisfy;
$$m_{i}<m_{1}+m_{2}+\cdots+\widehat{m_{i}}+\cdots+m_{n};\quad i=1,\ldots,n.$$\\
Under this condition variety of semi-stable configurations $\overline{\mathcal{C}}_{n}(m)$ is a projective compactification of $\mathcal{C}_{n}(m)$ by a finite number of points. Corresponding line bundle $\mathcal{L}$ of $\overline{\mathcal{C}}_{n}(m)$ can be written as 
$$\mathcal{L}(\Sigma)=\mathcal{L}_{1}^{\otimes m_{1}}\otimes \cdots \otimes \mathcal{L}_{n}^{\otimes m_{n}},$$\\
at a point $\Sigma=(p_{1},\ldots,p_{n}) \in \overline{\mathcal{C}}_{n}(m)$.\\
Setting $\mathcal{C}_{n}(m)=(\mathbb{P}_{1}^{n})^{s}/PSL_{2}(\mathbb{C}).$ We consider the map
$$\pi:(\mathbb{P}_{1}^{n})^{s}=\underbrace{(\mathbb{P}^{1}\times \cdots \times \mathbb{P}^{1})^{s}}_{n-copies} \longrightarrow \mathcal{C}_{n}(m).$$\\
With the fiber $\pi^{-1}(\Sigma)\simeq PSL_{2}(\mathbb{C})$. This is the structure of a principal $PSL_{2}(\mathbb{C})$ bundle.\\
Let $\xi$ be the linear vector bundle such that for $\Sigma \in \mathcal{C}_{n}(m),\ \xi(\Sigma)$ is the tangent space to fiber $\pi^{-1}(\Sigma)$ i.e, $\xi(\Sigma)$ is the tangent to $PSL_{2}(\mathbb{C})=\mathfrak{sl}_{2}$ which acts on $SL_{2}$ by $ad g:A\mapsto g^{-1}Ag$ and $PSL_{2}(\mathbb{C})=SL_{2}/ \pm 1.$\\
Note that $det\ ad g=1$ since $A\in SL_{2}$. So determinant bundle $det \xi$ is trivial. Taking into consideration all above, we form the Euler sequence to be
$$0\longrightarrow \xi \longrightarrow \bigoplus_{i=1}^{n}\mathcal{L}_{i}\longrightarrow \mathcal{T} \longrightarrow 0$$\\
where $\mathcal{T}$ is tangent bundle to $\mathcal{C}_{n}(m).$\\
The \textit{ canonical } bundle of $\mathcal{C}_{n}(m)$ is defined to be the determinant bundle of 1-forms $\Omega$ on $\mathcal{C}_{n}(m)$. Namely,
$$\kappa=det \Omega.$$\\
We know that $\Omega=\mathcal{T}^{*}$, dual of tangent bundle $\mathcal{T}$. Then we say $-\kappa=det \mathcal{T}$ is the \textit{anticanonical bundle}.\\
In an exact sequence 
$$0 \longrightarrow E^{\prime} \longrightarrow E \longrightarrow E^{\prime \prime} \longrightarrow 0$$
of vector bundles we have $detE=(detE^{\prime})\otimes(det E^{\prime \prime})$ and $det \bigoplus E_{i}=\bigotimes E_{i}$, again $E_{i}$`s are vector bundles. (see \cite{Ko} ). Then we arrive the following theorem \\
\textbf{Theorem :}  $-\kappa=det \mathcal{T}=\bigotimes_{i=1}^{n} \mathcal{L}_{i}$ where $\mathcal{T}$ is the tangent bundle of $\mathcal{C}_{n}(m)$ and $\mathcal{L}_{i}`s$ are natural bundles.\\
{\bf Proof:} By above argument and Euler sequence
$$0\longrightarrow \xi \longrightarrow \bigoplus_{i=1}^{n}\mathcal{L}_{i}\longrightarrow \mathcal{T}\longrightarrow 0$$\\
we have $det \bigoplus_{i=1}^{n}\mathcal{L}_{i}=det \mathcal{T}\otimes det \xi$ and $det \xi$ is trivial. So
$$det\mathcal{T}=\bigotimes_{i=1}^{n}\mathcal{L}_{i} \quad \Box$$.\\
\begin{defn} A topological space $X$ is called an even-cohomology space if its cohomology groups $H^{*}(X;\mathbb{Z})$ vanish for * odd.
\end{defn}
The following lemma is a first step to determine the cohomology of spatial polygons. For the proof see \cite{HK},\cite{KLY1}\\
\textbf{Lemma :} $\mathcal{M}_{n}$ is an even-cohomology space. \qquad $\Box$

As a consequence of the lemma, odd Betti numbers of $\mathcal{M}_{n}$ vanish. The following theorem is a useful tool for calculating Poincar\`{e} polynomials. Proof can be found in \cite{KLY1} and \cite{HK}\\
\textbf{Theorem :} Poincar\'e polynomial of the variety $\mathcal{M}_{n}$ is given by 
$$P_{q}(\mathcal{M}_{n})=\frac{1}{q(q-1)}((1+q)^{n-1}-\sum_{m_{I}\leq 
\frac{m}{2}} q^{|I|})$$\\
where $m=m_{1}+\cdots+m_{n}; m_{I}=\sum_{i\in I} m_{i}$.

Let us go back to variety of weighted stable configurations. For any decomposition $I\amalg J \amalg K \ldots =\{ 1,\ldots,n \}$ let $D_{IJK\ldots}$ be the cycle of stable configurations $\Sigma=(p_{1},\ldots,p_{n})$ with $p_{\alpha}=p_{\beta}$ for $\alpha, \beta$ are in the same component $I,J,K,\ldots$. In particular, we define
$$D_{ij}=\mbox{ divisor of stable configurations with } p_{i}=p_{j}.$$\\
We would like to characterize all effective cycles in $\mathcal{C}_{n}(m)$ using degenerate configurations $D_{IJKL \cdots}$.\\
\textbf{Theorem :} Any effective cycle in $\mathcal{C}_{n}(m)$ is equivalent to positive combinations of degenerate configurations $D_{IJKL \cdots}$.\\
{\bf Proof:} Theorem holds for special values for $m_{i}$`s. For example, for one massive point or three massive points. In these cases $\mathcal{C}_{n}(m) \simeq \mathbb{P}^{n-3}$ and $\mathcal{C}_{n}(m) \simeq (\mathbb{P}^{1})^{n-3}$ respectively.\\
It is possible to pass from one moduli space to another by a sequence of wall crossing $\mathcal{C}_{n}(m)\longrightarrow \mathcal{C}_{n}(\widetilde{m})$ such that only one inequality $m_{I}<\frac{m}{2}$ changes its direction to be $\widetilde{m_{I}}>\frac{m}{2}$ and all the other inequalities stay unchanged. In this case we may choose $m$ and $\widetilde{m}$ to be arbitrary close the wall $m_{I}=\frac{m}{2}$.\\
Let $\{ 1,\ldots,n \}=I\amalg J,\ I$ is the special subset mentioned above. Assume $|I|=k,\ |J|=l.$ Then
$$\mathcal{C}_{n}(m_{I},m_{i}:\ i\in J)\simeq \mathbb{P}^{l-2}\subset \mathcal{C}_{n}(m),$$
$$\mathcal{C}_{n}(\widetilde{m}_{J},m_{i}:\ i\in I)\simeq \mathbb{P}^{k-2}\subset \mathcal{C}_{n}(\widetilde{m})$$\\
and $\mathcal{C}_{n}(m)\slash \mathbb{P}^{l-2}$ is birationally equivalent to $\mathcal{C}_{n}(\widetilde{m})\slash \mathbb{P}^{k-2}$. Algebraic cycles in $\mathcal{C}_{n}(\widetilde{m})$ are those in $\mathcal{C}_{n}(m)$ and cycle in $\mathbb{P}^{l-2}$ are generated by degenerate configurations by the argument at the beginning of this proof.\\
$\Box$

Recall that $\mathcal{L}_{i}$`s are natural bundles on $\mathcal{M}$ such that the fiber at \linebreak $\Sigma=(p_{1},\ldots,p_{n})$ is the tangent space at $p_{i}\in \mathbb{P}^{1}$. Set
$$\l_{i}=[\mathcal{L}_{i}]=\{ \mbox{ zeros of }s \}-\{ \mbox{ poles of }s \}$$\\
where $s$ is a rational section of $\mathcal{L}_{i}$.\\
\textbf{Lemma :} With the previous notations $l_{i}=D_{ij}+D_{ik}-D_{jk}$ which is independent of choice of $j,k$.\\
{\bf Proof:} Let $t=\frac{p_{i}-p_{j}}{p_{i}-p_{k}}:\frac{z-p_{j}}{z-p_{k}}$ be local parameter at $z\in \mathbb{P}^{1}$ with $t(p_{i})=1$. Then
$$w_{i}=\frac{dt}{dz}=\frac{(p_{k}-p_{j})dz}{(z-p_{j})(z-p_{k})} \arrowvert_{z=p_{i}}=\frac{(p_{k}-p_{j})dp_{i}}{(p_{i}-p_{j})(p_{i}-p_{k})}$$\\
is the rational section of dual bundle $\mathcal{L}_{i}^{*}$. So $[w_{i}]=D_{jk}-D{ij}-D_{ik}$. Therefore,
$$l_{i}=D_{ij}+D_{ik}-D_{jk}.\quad \Box$$\\
\textbf{Corollary :} Some of the other relations between $l_{i}$ and $D_{ij}$ are as follows;
\begin{eqnarray}
1)&&D_{ij}=\frac{1}{2}(l_{i}+\l_{j}), \nn \\
2)&&l_{i}-l_{j}=\frac{1}{2}(D_{ik}-D_{jk}). \nn
\end{eqnarray}
{\bf Proof:} We have $l_{i}=D_{ij}+D_{ik}-D_{jk}.$ So \\
1)\begin{displaymath}
\left. \begin{array}{rr}
\mbox{ } & l_{i}=D_{ij}+D_{ik}-D_{jk} \\
\mbox{ } & l_{j}=D_{ij}+D_{jk}-D_{ik}
\end{array} \right \} \Rightarrow l_{i}+l_{j}=2D_{ij}
\end{displaymath}
so $D_{ij}=\frac{1}{2}(l_{i}+l_{j}).$\\
2) Follows from above.\qquad $\Box$\\

The lemma gives an inductive procedure to evaluate any monomial in $l_{i}$ in terms of "degenerate" cycles $D_{I,J,K,\ldots}$(in which all points $p_{i}\in I$ are glued together as well as for $J,K,L,\ldots$). The following corollary allows us to evaluate an arbitrary monomial in $l_{i}$. Note that non-zero cycles should contain at least three components and 3-component cycles represent a point provided $m_{I},m_{J},m_{K}$ satisfy triangle inequality.\\
\textbf{Corollary :}$l_{i}\cdot D_{I,J,K,\ldots}=D_{(IJ),K,\ldots}+D_{(IK),J,\ldots}-D_{I,(JK),\ldots}$.\\
{\bf Proof:} For any $l_{i}$ and cycle $D_{I,J,K,\ldots}$, with $i\in I$ we may write
$$l_{i}\cdot D_{I,J,K,\ldots}=[\mathcal{L}_{i}|\mathcal{C}_{n}(m_{I},m_{J},m_{K},\ldots)],$$\\
where right hand side of the above equation is the class of $\mathcal{L}_{i}$ in $\mathcal{C}_{n}(m_{I},m_{J},m_{K},\ldots) \simeq D_{I,J,K,\ldots}$ and $\mathcal{C}_{n}(m_{I},m_{J},m_{K},\ldots)$ is the moduli space of weighted stable configurations obtained from summing up weights of $\mathcal{C}_{n}(m)$ whose corresponding indices contained in $I,J,K,L,\ldots$\\
By lemma we may write
$$[\mathcal{L}_{i}|\mathcal{C}_{n}(m_{I},m_{J},m_{K},\ldots)]=D_{IJ}+D_{IK}-D_{JK}$$\\
provided $i \in I$ and we obtain the equality in the statement of corollary.\quad $\Box$\\
\textbf{Example :} For $i\neq j$ we can evaluate $l_{j}l_{i}$ as follows; we know that $l_{i}=D_{ij}+D_{ik}-D_{jk}.$ So,
\begin{displaymath}
\left. \begin{array}{rr}
\mbox{ } & l_{j}\cdot D_{ij}=D_{ijk}+D_{ijl}-D_{(ij)(kl)} \\
\mbox{ } & l_{j}\cdot D_{ik}=D_{ijk}+D_{(ik)(jl)}-D_{ikl} \\
\mbox{ } & l_{j}\cdot D_{jk}=D_{ijk}+D_{jkl}-D_{(il)(jk)}
\end{array} \right.
\end{displaymath}
This implies
$$l_{j}l_{i}=D_{ijk}+D_{ijl}-D_{ikl}-D_{jkl}+D_{(ik)(jl)}+D_{(il)(jk)}-D_{(ij)(kl)}.$$\\
\textbf{Example :}For $p=l_{i}^{2}$, a similar calculation leads us to the formula
$$p=D_{ijk}+D_{ijl}+D_{ikl}+D_{jkl}-D_{(ij)(kl)}-D_{(ik)(jl)}-D_{(jk)(il)}.$$\\
This expression is independent of $i,j,k,l.$\\

By the equivalence of stable configurations and spatial polygons, we can relate the divisors $D_{ij}$ by some kind of polygons. In other words, a divisor $D_{ij}$ corresponds to a polgon in $\mathcal{M}$ with edges $(p_{1},p_{2},\ldots,p_{n})$ and $p_{i}\uparrow \uparrow p_{j}$ i.e, $p_{i}$ and $p_{j}$ are parallel. For anti-parallel edges we write $p_{i} \uparrow \downarrow p_{j}$.\\
Using the following theorem we may calculate the cohomology rings of stable configurations, hence cohomology rings of spatial polygons.
\textbf{Theorem :} The Chow(cohomology) ring of $\mathcal{C}_{n}(m)$ is generated by the class of divisors $D_{ij}$ subject to the following relations;\\
1) $\forall$ quadruple $(i,j,k,l)$ there are linear relations
$$D_{ij}+D_{km}=D_{ik}+D_{jm}=D_{im}+D_{kj}=\frac{1}{2}(l_{i}+l_{j}+l_{k}+l_{m}).$$\\
2) For any triple $(i,j,k)$ there are quadratic relations
$$D_{ij}D_{jk}=D_{jk}D_{ki}=D_{ki}D_{ij}=D_{ijk}.$$\\
3) For any tree $\Gamma$ with vertices in $I\subset \{ 1,\ldots,n \}$ such that $m_{I}>\frac{m}{2},$
$$D_{I}=\prod_{(ij)\in \Gamma} D_{ij}=0.$$\\
{\bf Proof:} We know that the divisor $D_{ij}$ generate the Chow ring. In the view of the formula $D_{ij}=\frac{1}{2}(l_{i}+l_{j})$, relations in 1) becomes trivial. Using the same formula we also see that all products in quadratic relations 2) are equal to $D_{ijk}$.\\
The product in 3) is a locus of configurations with equal points $p_{i},i\in I$. Under the condition $m_{I}>\frac{m}{2}$, such configuration is unstable and hence $\prod_{(ij)\in \Gamma} D_{ij}=0$. Observe that the quadratic relations ensure that the product $D_{I}$ is independent of choice of tree $\Gamma$ on vertices $I$.\\
To prove the completeness, we need to show that; for any disjoint subsets $I,J,K,M \subset \{ 1,\ldots,n \}$ we have
$$(*)\quad D_{(IJ),K,M}+D_{I,J,(KM)}=D_{(IK),J,M}+D_{I,K,(JM)}=D_{(IM),J,K}+D_{I,M,(JK)}.$$\\
In fact, if $i,j,k,m$ are elements from $I,J,K,M$ respectively then the above equation is equivalent to the following identities;
$$D_{I,J,K,M}(D_{ij}+D_{km})=D_{I,J,K,M}(D_{ik}+D_{jm})=D_{I,J,K,M}(D_{im}+D_{jk})$$\\
which follows from 1). So (*) holds.\\
Now, let us consider a puzzle; let`s divide a heap of stones of masses $m_{i}$ into three parts of masses $m_{I},m_{J},m_{K}$ satisfying the triangle inequality. Then any other such division may be obtained from the initial one br removing a stone from one heap and putting it into another so that new heaps also satisfy the triangle inequality.\\
Using the puzzle we can show that if $I\amalg J\amalg K=\{ 1,\ldots,n \}$ is stable decomposition i.e, $m_{I},m_{J},m_{K}$ satisfy the triangle inequality. Then $D_{I}D_{J}D_{K}=D_{IJK}$ is independent of stable decomposition.\\
Really, by the puzzle it is enough to check that
$$D_{I\slash \{ i \},J\cup \{ i \},K}=D_{IJK} \mbox{ if } I\slash \{ i \},J\cup \{ i \},K \mbox{ is stable.}$$\\
Applying (*) to the quadruple $\{ i \},I\slash \{ i \},J,K$ we get
$$D_{IJK}+D_{\{ i \},I\slash \{ i \},(JK)}=D_{(\{ i \}J),K,I\slash \{ i \}}+D_{\{ i \},J,(KI\slash \{ i \})}.$$\\
Using triangle inequalities $m_{J}+m_{K}>\frac{m}{2}$ and $m_{K}+M_{I\slash \{ i \}}>\frac{m}{2}$ we obtain unstable decompositions $(JK)$ and $(KI\slash \{ i \})$ then
$$D_{\{ i \},I\slash \{ i \},(JK)}=D_{\{ i \},J,(KI\slash \{ i \})}=0$$.\\
 So we obtained the desired result.
It remains to check that\\
\textit{i) }Linear relations between divisors are complete,\\
\textit{ii) }Linear relation on divisor  cycles $D_{IJK\ldots}$ follows from the relation 1) in the statement of the theorem.\\
The second part \textit{ii) } is the consent of (*). The first part \textit{i) } necessarily says that the cross ratio $[p_{i}:p_{j}:p_{k}:p_{l}]$ generate the whole ring of non-vanishing regular functions on the divisors $D_{ij}$. Actually, if we fix $p_{i}=0, p_{j}=\infty, p_{k}=1$ then the complement of the divisor $D_{ij}$ is a subset $X\subset \mathbb{C}^{n-3}$ with pairwise distinct components$\neq 0,1$. Non-vanishing regular functions on $X$ are generated by $(p_{l}-p_{m})^{\pm 1},p_{l}^{\pm 1},(p_{l}-1)^{\pm 1}$ and may be expressed by cross-ratio. Hence we are done.\quad $\Box$\\
\textbf{Corollary :}The Chow(cohomology) ring $H^{*}(\mathcal{M})$ over $\mathbb{Z}$ is generated by the classes of natural bundles subject to relations
\begin{eqnarray}
1)&& l_{i}^{2}=p,\mbox{ independent of } i \nn \\
2)&& \sum_{2k+r=|I|-1}p^{k}\sigma_{r}(l_{I})=0,\quad m_{I}>\frac{m}{2},\nn 
\end{eqnarray}
where $\sigma_{r}$ is the $r$-th elementary symmetric polynomial and $m=m_{1}+\cdots +m_{n};m_{I}=\sum_{i\in I} m_{i}.$

Above corollary gives a handy method to determine the Chow ring of polygon spaces. Moreover, in what follows we will give explicit formulae to determine the monomials of the form $p^{k}l_{J},J\subset I$ with $m_{I}>\frac{m}{2}.$

{\bf Non-Holomorphic Cycles.}
There are natural symplectic cycles in the polygon space $\mathcal{M}$, corresponding degenerate polygons with antiparallel sides. Tangent space to the moduli of such degenerate polygons is not closed under complex multiplication but as a topological cycle, it should be an integer linear combination of the holomorphic cycles $D_{I,J,K,\ldots}.$ In a simplest case of non-holomorphic divisor, it may be expressed as follows,\\
\textbf{Lemma :} The cycle $D^{-}_{ij}$(see Fig.9b) of polygons with antiparallel vectors $p_{i}\uparrow \downarrow p_{j}$ is equivalent to $\frac{1}{2}(l_{i}-l_{j})$ if $D^{-}_{ij}$ is oriented by the vector $p_{j}$.\\
{\bf Proof:} Let us compare the intersection indices of $D^{-}_{ij}$ and $\frac{1}{2}(l_{i}-l_{j})$ with holomorphic curves, that is with quadrangles $D_{IJKL}$. Both intersections $D^{-}_{ij}\cdot D_{IJKL}$ and $\frac{1}{2}(l_{i}-l_{j})\cdot D_{IJKL}$ are zero if $i,j$ belongs to the same set $I,J,K,L$. So we may suppose that $i\in I,j\in J$. Then the first intersection index is non-zero if and only if the quadrangle admits a degeneration into triangle with antiparallel vectors $p_{I}\uparrow \downarrow p_{J}$. Here $p_{I}=\sum_{i\in I} p_{i}.$ In this case
$$||p_{K}-p_{L}||<||p_{I}-p_{J}||<p_{K}+p_{L} \mbox{ i.e, }|m_{K}-m_{L}|<|m_{I}-m_{J}|<m_{L}+m_{K}.$$\\
Assume $m_{I}>m_{J}$ and $m_{K}>m_{L}$. Then previous inequality becomes
$$m_{K}+m_{J}<m_{I}+m_{L},m_{I}<m_{J}+m_{K}+m_{L}. \quad (*)$$\\
In the previous section we saw that there are two types of quadrangles(shown in Fig.1) "triangle" type and "star" type which are represented by $(ijk)$ and $\{ (ij)(ik)(il) \}$ respectively. In our case we have 
$$m_{I}+m_{K}>m_{J}+m_{L} \mbox{ and }m_{K}+m_{J}<m_{I}+m_{L}.$$\\
 So if $m_{I}+m_{J}>m_{K}+m_{L}$ we obtain the star $\{ (IJ)(IK)(IL) \}$, otherwise we find the triangle $(IKL)$. By corollary(1.40) we have 
$$l_{i}\cdot D_{IJKL}=D_{(IJ)KL}+D_{(IK)JL}-D_{I(JK)L}$$.\\
Together with this formula and using the assumption $i\in I, j\in J$ we have\\
\begin{eqnarray}
l_{i}\cdot D_{IJKL}=0,&& l_{j}\cdot D_{IJKL}=2 \ \mbox{ for triangle }(IKL),\nn \\
l_{i}\cdot D_{IJKL}=-1,&& l_{j}\cdot D_{IJKL}=1 \ \mbox{ for star }\{ (IJ)(IK)(IL) \}. \nn
\end{eqnarray}
In both cases $\frac{1}{2}(l_{i}+l_{j})\cdot D_{IJKL}=-1$.\\
Suppose $m_{I}<m_{J}$ then intersection indices change sign. As a result we get,
\begin{displaymath}
\frac{1}{2}(l_{i}-l_{j})\cdot D_{IJKL}=
\left \{ \begin{array}{ll}
1 & \mbox{ if } m_{I}<m_{J}\\
-1 & \mbox{ if } m_{I}>m_{J}.
\end{array} \right.
\end{displaymath}
To compare this to $D^{-}_{ij}\cdot D_{IJKL}$ we must choose an orientation of $D^{-}_{ij}$, which is in fact an orientation of the normal plane to $p_{i}$ or $p_{j}$. Let us choose the orientation by the normal to $p_{i}$. Then 
\begin{displaymath}
D^{-}_{ij}\cdot D_{IJKL}=
\left \{ \begin{array}{ll}
1 & \mbox{ if } m_{I}>m_{J}\\
-1 & \mbox{ if } m_{I}<m_{J}
\end{array} \right. =-\frac{1}{2}(l_{i}-l_{j})\cdot D_{IJKL}.
\end{displaymath}
So if we orient it by $p_{j}$ we obtain the Lemma.\quad $\Box$

This lemma is quite useful in determining any product
$$D_{I}^{\epsilon}=2^{1-|I|}\prod_{(ij)\in \Gamma} (l_{i}\pm l_{j})$$\\
where $\epsilon(i,j)=-\epsilon(j,i)=\pm 1$ and $m_{I}>\frac{m}{2}$ with a tree $\Gamma$ with vertices $I$. In a different point of view, we may consider $\Gamma$ as partially oriented i.e, some edges have an orientation $i\rightarrow j$, yet the others may not be oriented. So we define
$$D_{\Gamma}=2^{1-|I|}\prod_{(ij)\in \Gamma} (l_{i}\pm l_{j})=2^{1-|I|} \prod (l_{i}+l_{j})\prod_{i\rightarrow j}(l_{i}-l_{j}).$$\\
Let us take $\Gamma$ to be star $\{ (i_{0}i_{1})(i_{0}i_{2})\cdots (i_{0}i_{q}) \}$. Then
$$D_{I}^{\epsilon}=2^{1-|I|}\prod_{i\in I\slash \{ i_{0} \} }(l_{i_{0}}+\epsilon_{i}l_{i})$$\\
with $\epsilon_{i}=\pm 1$ and $\epsilon_{i_{0}}=1.$\\
Observe that
$$\prod_{(ij)\in \Gamma} (l_{i}+l_{j})=\sum_{\tiny \begin{array}{c}
orientations\ of \\
\Gamma \end{array}}\prod_{i\in I} l_{i}^{deg^{+}i}=\sum_{2k+|I_{odd}|=|I|-1} p^{k}l_{I_{odd}},$$\\
where $deg^{+}i$ is the number of edges entering $i$ and $I_{odd}$ is a subset of vectors with odd $deg^{+}$. This leads to the equation
$$D_{I}^{\epsilon}=2^{1-|I|} \sum_{\tiny \begin{array}{c}
J\subset I\\
|I\slash J|=2k+1 \end{array}} \epsilon_{J}l_{J}p^{k} \quad (*)$$\\
where $\epsilon_{J}=\prod_{i\in J}\epsilon_{i}, l_{j}=\prod_{j\in J}l_{j}.$\\
We may take the last formula as a system of $2^{1-|I|}$ equations with the same number of variables $l_{J}p^{k}, J\subset I$ and $|I\slash J|=2k+1$. Note that the matrix of the system is invertible because its square is a scalar matrix.\\
The following theorem and its corollary are useful in calculating the cohomology.
\textbf{Theorem :} For any $J\subset I,m_{I}>\frac{m}{2}$ and $|I\slash J|=2k+1$, the following formula holds
$$l_{J}\cdot p^{k}=\sum_{\epsilon} \epsilon_{J}D_{I}^{\epsilon}$$\\
where the product is taken over all combinations of signs $\epsilon_{i}=\pm1,\epsilon_{0}=1.$\\
{\bf Proof:} Using (*) we can write
$$\sum_{\epsilon} \epsilon_{J}D_{I}^{\epsilon}=2^{1-|I|}\sum_{\tiny \begin{array}{c}
K\subset I\\
|I\slash K|=2k+1 \end{array}}(\sum_{\epsilon}\epsilon_{J}\epsilon_{K})l_{K}p^{k}=l_{J}p^{k}$$\\
Since 
\begin{displaymath}
\sum_{\epsilon}\epsilon_{J}\epsilon_{K}=\left \{ \begin{array}{ll}
0 & \mbox{ if }J\neq K, \\
2^{|I|-1} & \mbox{ if }J=K .
\end{array} \right. \qquad \Box
\end{displaymath}
\textbf{Corollary :}Let $|J|+2k=n-3=dim \mathcal{M}$ and $I\supset J$ be any set of cardinality $n-2$, say $I=\{ 1,\ldots,n \}\slash \{ \alpha,\beta \},\alpha,\beta \notin I.$ Then 
$$l_{J}p^{k}=\sum_{|m_{\alpha}-m_{\beta}|<|(\epsilon_{I},m_{I})|<m_{\alpha}+m_{\beta}} sgn(\epsilon_{I},m_{I})\epsilon_{I\slash J}$$\\
where the sum is taken over all signs $\epsilon_{i}=\pm 1$, fixed on one element $\gamma \in I$ with $\epsilon_{\gamma}=1$ and $(\epsilon_{I},m_{I})=\sum_{i\in I}\epsilon_{i}m_{i}.$\\
{\bf Proof:} By the assumption of the corollary $|I|=n-2$ and hence $|I\slash J|=2k+1$. Note that in the sum, we consider $\alpha,\beta \in \{1,\ldots,n \} \slash I$ such that the polygon degenerates into a triangle and satisfies the triangle inequality. So the previous theorem is applicable. $\quad \Box$
\section{Fano Polygon Spaces}
\textbf{Proposition :}Let the vertices of an element of the moduli space $\mathcal{M}$ be numbered as $1,2,\ldots,n$. Then the first Chern class $c_{1}(\mathcal{M})$ is given by 
$$c_{1}(\mathcal{M})=\sum_{i=1}^{n} D_{ii+1}.$$\\
Here we set $D_{nn+1}=D_{n1}$\\
{\bf Proof:} Using the natural bundles we write $c_{1}(\mathcal{M})=\sum_{i}l_{i}$. Taking consecutive sum with the convention $l_{n+1}=l_{1}$ we have,
$$c_{1}(\mathcal{M})=\sum_{i}l_{i}=\sum_{i=1}^{n}\frac{1}{2}(l_{i}+l_{i+1})=\sum_{i=1}^{n}D_{ii+1}.\quad \Box$$\\
\begin{defn} A divisor $D$ is ample if $D\cdot C>0,\ \forall$ curve $C$ contained in $\mathcal{M}_{n}$. Similarly, a vector bundle $\mathcal{L}$ is ample if $\mathcal{L}|_{C}$ is ample for every curve in $\mathcal{M}_{n}.$
\end{defn}
Note that $D_{IJKL}\neq 0 \Leftrightarrow m_{I},m_{J},m_{K},m_{L}$ satisfy $m_{I}<m_{J}+m_{K}+m_{L}$ for all $I,J,K,L$ and in this case $D_{IJKL} \simeq \mathbb{P}^{1}.$\\
\textbf{Theorem ;}[Ampleness Criterion] $D=\sum_{i=1}^{n} a_{i}l_{i}$, where $l_{i}$`s are characteristic classes of natural bundles $\mathcal{L}_{i}$ and $a_{i}$`s are positive numbers, is ample if and only if for all quadruple $D_{IJKL}$ we have
\begin{displaymath}
\left \{ \begin{array}{ll}
a_{I}>0 & \mbox{ for "triangle" }\\
a_{I}<a_{J}+a_{K}+a_{L} & \mbox{ for "star"}
\end{array} \right.
\end{displaymath}
{\bf Proof:} Recall that there are two types of quadrangles;\\
1) Triangle type; for indices $i,j,k,l$ 
$$m_{i}+m_{j}>\frac{m}{2},\ m_{j}+m_{k}>\frac{m}{2},\ m_{k}+m_{i}>\frac{m}{2}$$\\
2) Star type; for indices $i,j,k,l$\\
$$m_{i}+m_{j}>\frac{m}{2},\ m_{i}+m_{k}>\frac{m}{2},\ m_{i}+m_{l}>\frac{m}{2}$$\\
For triangle type, $l_{i}\cdot D_{ijkl}=D_{ij}+D_{ik}-D_{jk}=0$ and $l_{l}\cdot D_{ijkl}=D_{li}+D_{lj}-D_{ij}=2$\\
For star type ,$l_{i}\cdot D_{ijkl}=D_{ij}+D_{ik}-D_{jk}=-1$ and $l_{j}\cdot D_{ijkl}=D_{ji}+D_{jk}-D_{ik}=0.$\\
We need to check $D\cdot D_{IJKL}=\sum_{i=1}^{n} a_{i}(l_{i}\cdot D_{IJKL})>0,\ \forall D_{IJKL}\neq 0.$ So we will have two cases,\\
\textbf{Case I} Triangle case: In this case we set
$$m_{L}+m_{J}>\frac{m}{2},\ m_{L}+m_{K}>\frac{m}{2},\ m_{J}+m_{K}>\frac{m}{2}$$\\
and
\begin{displaymath}
l_{i}\cdot D_{IJKL}=
\left \{ \begin{array}{ll}
2 & \mbox{ for } i \in I \\
0 & \mbox{ for } i\notin I
\end{array} \right.
\end{displaymath}
\textbf{Case II} Star case: In this case
$$m_{I}+m_{L}>\frac{m}{2},\ m_{I}+m_{J}>\frac{m}{2},\ m_{I}+m_{K}>\frac{m}{2}.$$\\
So
\begin{displaymath}
l_{i}\cdot D_{IJKL}=
\left \{ \begin{array}{ll}
-1 & \mbox{ for } i \in I \\
1 & \mbox{ for } i\notin I
\end{array} \right.
\end{displaymath}
In order to obtain $D\cdot D_{IJKL}>0$, in the first case we must have $a_{I}>0$ and in the second case $a_{I}<a_{J}+a_{K}+a_{L}. \Box$\\
\textbf{Corollary :}
The anticanonical class $c_{1}(\mathcal{M}_{n})=\sum l_{i}$ is ample if and only if for all star type degenerations of any quadruple $D_{IJKL}$ we have
$$|I|<|J|+|K|+|L|.$$\\
\begin{defn} $\mathcal{M}$ is Fano if the first anti-canonical class(first Chern class) is ample.
\end{defn}

In other words, $\mathcal{M}$ is Fano if and only if $c_{1}(\mathcal{M})\cdot D_{IJKL}>0$ for all quadrangles $D_{IJKL}$.\\
\begin{defn} A maximal degeneration in $\mathcal{C}_{n}(m)$ is a cycle consisting of configurations in which $p_{i}=p_{j}$, for all $i,j\in I$ and $I$ is the maximal set. We will denote a maximal degeneration by $\mathcal{M}_{I}$.
\end{defn}

Note that $\mathcal{M}_{I}\subset \mathcal{C}_{n}(m)$ and actually $\mathcal{M}_{I} \simeq \mathbb{P}^{k}$ with $k=n-|I|-2$.\\
The maximal sets $I$ are characterized by inequalities 
\begin{eqnarray}
m_{I}+m_{s}&>&\frac{m}{2} \nn \\
m_{I}&<&\frac{m}{2} \nn
\end{eqnarray}\\
\textbf{Theorem :}$\mathcal{C}_{n}(m)$ is Fano if and only if any maximal degeneration $\mathcal{M}_{I}$ is either a point or has dimension greater than $\frac{n-4}{2}$\\
{\bf Proof:} We know that $\mathcal{C}_{n}(m)$ is Fano $\Leftrightarrow$ the polygon has no quadrangle degenerations of "star" type i.e, we won`t have degenerations given  by 
$$m_{I}+m_{J}>\frac{m}{2},\ m_{I}+m_{K}>\frac{m}{2},\ m_{I}+m_{L}>\frac{m}{2}$$\\
with $|I|\geq \frac{n}{2}.$\\
We can successively move a side $p_{s}, \ s\notin I$ to the set of edges whose indices are contained in $I$ as long as it is possible i.e, we can move it as long as $m_{I}<\frac{m}{2}$. As a result, we arrive to the maximal degeneration $m_{I}<\frac{m}{2}$ and $m_{I}+m_{s}>\frac{m}{2},\ \forall s\notin I.$ There are two possibilities;\\
\textit{i)} The maximal degeneration is a point i.e, $|I|=n-2$,\\
\textit{ii)} The maximal degeneration is of positive dimension.
$$dim \mathcal{M}_{I}=n-|I|-2<\frac{n}{2}-2$$\\
for the ones which are not Fano(i.e, quadrangles of star type with $|I| \geq \frac{n}{2})$\\
Therefore, for Fano we have 
$$dim \mathcal{M}_{I}>\frac{n-4}{2}.\quad \Box$$\\
\textbf{Corollary :} For $n\in \{ 4,5 \}$ we have Fano polygon spaces\\

Vehbi Emrah Paksoy, Claremont McKenna College, emrah.paksoy@cmc.edu

\end{document}